\documentclass{article}
\usepackage{hyperref}
\hypersetup{
  colorlinks   = true, 
  urlcolor     = black, 
  linkcolor    = blue, 
  citecolor   = green 
}
\usepackage{geometry}
\usepackage[english]{babel}
\usepackage{amsmath}
\usepackage{graphicx}
\usepackage{ marvosym }
\usepackage[utf8]{inputenc}
\usepackage{mathtools}
\usepackage{geometry}
\usepackage{amsfonts}
\usepackage{amssymb}
\usepackage{amsthm}
\usepackage{thmtools}
\usepackage{t1enc}
\usepackage[titles]{tocloft}
\usepackage{makeidx}
\usepackage{wasysym}
\usepackage{stmaryrd}
\usepackage{calc}  
\usepackage{enumitem} 
\usepackage[refpage]{nomencl}

\usepackage{bookmark}
\usepackage{algpseudocode}
\usepackage{tikz}
\usetikzlibrary{arrows}
\usetikzlibrary{decorations}
\usepackage{float}
\usepackage{mathtools}

\usepackage{enumitem}
\usepackage{linegoal}
\usepackage{calc}

\theoremstyle{plain}
\newtheorem{thm}{Theorem}

\newtheorem{claim}[thm]{Claim}

\newtheorem{prop}[thm]{Proposition}
\newtheorem*{prop*}{Proposition}
\newtheorem*{seged*}{Sublemma}
\newtheorem{cor}[thm]{Corollary}
\newtheorem{lem}[thm]{Lemma}

\newtheorem*{lem*}{Lemma}
\theoremstyle{definition}

\newtheorem*{defn*}{Definition}

\newtheorem{fel*}[thm]{Exercise}

\newtheorem*{megf*}{Observation}
\theoremstyle{remark}
\newtheorem{rem}[thm]{Remark}
\newtheorem*{rem*}{Remark}

\newenvironment{sbiz}{\par\noindent{\itshape Proof:}\ }{\rule{1ex}{1ex}}
\newenvironment{ssbiz}{\par\noindent{\itshape Proof:}\ }{\rule{1ex}{1ex}\ \rule{1ex}{1ex}}

\newenvironment{nbiz}{\par\noindent{\itshape Proof:}\ }{}

\title{Gomory-Hu trees of infinite graphs with finite total weight }
\author{Attila Joó \thanks{MTA-ELTE Egerváry Research Group, Department of Operations Research, Eötvös Loránd University,
Budapest, Hungary.
 Email: {\tt joapaat@cs.elte.hu} }}
\date{2016}

\begin{document}
\maketitle

\begin{abstract}
 Gomory and Hu proved  in \cite{gomory1961multi} their well-known theorem which states that if $ G $ is a finite graph with non-negative 
 weights on its edges, then there exists a tree $ T $ (called now Gomory-Hu tree) on $ V(G) $ such that for all $ u\neq v\in V(G) $ there 
 is an $ e\in E(T) $ such that the two components of $ T-e $ determines an optimal (minimal valued) cut between $ u $ an $ v $ in $ G 
 $. In this paper we extend their result to infinite weighted graphs with finite total 
  weight. Furthermore, we show by an example that one can not omit the condition of finiteness of the total weight.    
\end{abstract}

\section{Introduction}

Let $ G=(V,E) $  be a countable connected simple graph and let $ 
c:E\rightarrow\mathbb{R}_{+}\setminus \{ 0 \} $ be a 
weight-function, then $ (V,E,c) $ is  a \textbf{weighted graph}. 
We call the subsets $ X $ of $ V $ cuts and we write $ \boldsymbol{\mathsf{out}_G(X)} $ for the set of the edges with exactly one end in 
$ X $. We say 
$ X $ is an $ \boldsymbol{u-v} $ \textbf{cut} for some 
$ u\neq v\in V $ if $ u\in X $ and $ v\notin X 
$. A cut $ X $ separates $ u $ and $ v $ if $ X $ is either a $ u-v $ or a $ v-u $ cut.  Let
 $ \boldsymbol{d_c(X)}=\sum_{e\in \mathsf{out}_G(X)} c(e)$ and  let   $ \boldsymbol{\lambda_c(u,v)}:=  \inf \{ d_c(X):\ 
X\text{ is a } 
u-v\text{ cut } \} $ for $ u\neq v\in V $. A cut $ X $ 
is an \textbf{optimal} $ \boldsymbol{u-v} $ \textbf{cut} if it is a $ 
u-v $ cut with $ d_c(X)=\lambda_c(u,v) $. A cut $ X $ is \textbf{optimal} if it is an optimal $ u-v $ cut for some $ u\neq v\in V $. The 
weighted graph $ (V,E,c) $ is  \textbf{finitely separable} 
if   $ \lambda_c $ has just finite values.  A tree $ T=(V,F) $ is a \textbf{Gomory-Hu tree} for $ (V,E,c) $ if for all 
 $ u\neq v\in V $ there is an $ e\in F $ such that the \textbf{fundamental cuts corresponding to} $ \boldsymbol{e} $ (i.e. the vertex 
 sets 
 of the components of $ 
 T-e $) separate optimally $ u $ and $ v $ in $ (V,E,c) $. Gomory and Hu proved in \cite{gomory1961multi} that for all finite weighted 
 graph  
 there exists a Gomory-Hu tree. It has several interesting consequences. For example the function $ \lambda_c $ may have at most $ n-1 $ 
 different  
 values instead of $ \binom{n}{2} $ (where $ n $ is the number of the vertices) and there is at least two optimal cuts that consist of a 
 single 
 vertex, namely the leafs of the  Gomory-Hu tree  (unless the graph is trivial).  
 
 In this paper we 
  extend their theorem for infinite weighted graphs with finite total weight. Note that, the strict positivity of $ c $ and the 
  connectedness 
  of $ G $  are not  real restrictions  since 
    throwing away edges $ e $ with $ c(e)=0 $ has no effect on the values of the cuts and one can construct Gomory-Hu trees 
    component-wise and join them to a Gomory-Hu tree. Furthermore, if the sum of the weights is finite, then the weighted graph must be 
    countable. 
    
    The cut structure of infinite graphs has been already investigated in some other perspectives (see for example 
    \cite{evangelidou2014cactus} and \cite{dunwoody2013structure}) where the authors only allow cuts with finitely many outgoing edges. 
    As it seems from the definitions above we are focusing on the literal generalization of Gomory-Hu trees. 
  In a more abstract folklore  version of the Gomory-Hu  theorem there is no 
  graph, one just has a finite set $ V $ and a  function $ 
     b:\mathcal{P}(V)\rightarrow \mathbb{R}_+ $ which is symmetric ($ b(X)=b(V\setminus 
      X) $) and submodular i.e.   \[ b(X)+b(Y)\geq b(X\cap Y)+b(X\cup Y)  \text{ if }  X, Y\subseteq V. \] Let $\lambda_b(u,v)= 
      \inf\{b(X): X \text{ a } u-v  \text{ cut}  \} $ (this definition makes sense for infinite $ V $ as well). In this case, there  
      exists an 
     \textbf{abstract Gomory-Hu tree}  with respect  to $ b $ in the following sense. There is a tree $ T $ on the vertex set $ V $ in 
     such a way that for 
     every  $ u\neq v\in V $ there is some 
          $ e\in E(T) $ such that for a fundamental cut $ X $ corresponding to $ e $, we have $ b(X)=\lambda_b(u,v). $

 \section{Preparations}
 Let $ (V,E,c) $ be a weighted graph.
 
  \begin{prop}\label{d_c szubmod}
   $ d_c(X)+d_c(Y)\geq d_c(X\cup Y)+d_c(X\cap Y) $ for all $ X,Y\subseteq V $.   
  \end{prop}
 \begin{sbiz}
   If edge $ e $ goes between $ X\setminus Y $ and $ Y\setminus X $, then it contributes  $ 2c(e) $ to the left side and  $ 
   0 $ to the right side of the inequality. The contribution of any other type of edge is the same for both sides. 
  \end{sbiz}
  
  For a sequence $ (X_n) $ let
  
     \begin{align*}
  &\liminf X_n= \bigcup_{m=0}^{\infty} \bigcap_{n=m}^{\infty} X_n\\
  &\limsup X_n= \bigcap_{m=0}^{\infty} \bigcup_{n=m}^{\infty} X_n.
  \end{align*}
  
  If $\liminf X_n= \limsup X_n $, then we denote this set by $ \lim X_n $ and we say that $ (X_n) $ is convergent.  
  \begin{claim}\label{limeszhalmaz költsége}
   
   \begin{minipage}[t]{\linegoal}
   \begin{enumerate}[leftmargin=*]
   \item If $ (X_n) $ is a convergent sequence of cuts, then $d_c(\lim 
    X_n)\leq \liminf d_c(X_n) $.
   \item In addition, if $ \sum_{e\in E}c(e)<\infty $, then $  \lim d_c(X_n)$  exists 
     and $ \lim d_c(X_n)=d_c(\lim X_n) $ holds.
   \end{enumerate}
   \end{minipage}
    
   \end{claim}
   
   \begin{nbiz}
   It is routine to check that $  \mathsf{out}_G(\lim X_n)= \lim \mathsf{out}_G(X_n)  $ holds. Consider the discrete measure space $ 
    (E,\mathcal{P}(E), \tilde{c}) $ where $ \tilde{c}(F)=\sum_{e\in F}c(e) $ for $ 
     F\subseteq E $.
     
      By applying the Fatou lemma to the characteristic functions of the sets $ 
    \mathsf{out}_G(X_n) $  we obtain
    
    \begin{align*}
    & d_c(\lim (X_n))=\widetilde{c}(\mathsf{out}_G(\lim X_n) )=\widetilde{c}(\lim \mathsf{out}_G(X_n))=\\ 
    &\widetilde{c}(\liminf \mathsf{out}_G(X_n)) \leq\liminf 
    \widetilde{c}(\mathsf{out}_G(X_n))= 
    \liminf d_c(X_n). 
    \end{align*} 
    
     At the  statement 2 we have $ \widetilde{c}(E)<\infty $,  thus the constant $ 1 $ function is integrable. Therefore by using 
     Lebesgue 
     theorem to 
     the 
     characteristic functions of the sets $ 
      \mathsf{out}_G(X_n) $ we obtain
      
      \[ d_c(\lim X_n)=\widetilde{c}(\lim \mathsf{out}_G(X_n))=\lim \widetilde{c}(\mathsf{out}_G(X_n))= \lim d_c(X_n).\ \rule{1ex}{1ex} \]
   \end{nbiz}
  Let us formulate our main result in the following  more abstract way.    
  \begin{thm}\label{GH van absztract GHtree}
      Let $ V $ be a nonempty countable set and let $ b: 
        \mathcal{P}(V)\rightarrow 
        \mathbb{R}_+\cup \{ \infty \} $ such that
       
       \begin{enumerate}
       \item[0.] $ b(X)=0 \Longleftrightarrow X\in \{ \varnothing,V \} $,
       \item[1.] $ b(X)=b(V\setminus X) \text{ for } X\subseteq V, $ ($ b $ is symmetric) 
       \item[2.] $ b(X)+b(Y)\geq b(X\cap Y)+b(X\cup Y)  \text{ for }  X, Y \subseteq V, $ ($ b $ is submodular)
       \item[3.]  if $ (X_n) $ is a nested sequence of cuts, then  $ b(\lim X_{n})=\lim b(X_n) $ ($ b $ is monotone-continuous)
       \item[4.]  $ \lambda_b $ has only finite values. ($ b $ finitely separate)
       \end{enumerate}
       Then there exists an abstract Gomory-Hu tree with respect to $ b $.
   \end{thm}
   
   \begin{rem}\label{GH szubmodegye 2}
   Properties 1,2 imply that for any $ X,Y $  we also have  \[ b(X)+b(Y)\geq b(X\setminus Y)+b(Y\setminus X). \]
   \end{rem}

 Observe that if  $ \sum_{e\in E}c(e)<\infty $ holds, then  $ b:=d_c $ 
 satisfies the properties above. Hence as a special case of Theorem \ref{GH van absztract GHtree} we obtain:
 \begin{cor}
 Every   weighted graph  with $ \sum_{e\in E}c(e)<\infty $ admits a Gomory-Hu tree. 
 \end{cor}
 
Consider the following weakening of 3.

\begin{enumerate}
\item[3']  if $ (X_n) $ is a nested sequence of cuts, 
   then  $ b(\lim 
   X_{n})\leq\liminf b(X_n) $.
\end{enumerate}
If we do not assume $ \sum_{e\in E}c(e)<\infty $ and we demand just $ (V,E,c) $ to be finitely separable, then Claim \ref{limeszhalmaz 
költsége} ensures 
that  $ b:=d_c $ still satisfy  this weaker condition (see Claim \ref{limeszhalmaz 
költsége}/1). 
We will see by a counterexample that in this case one can not guarantee the existence of a Gomory-Hu tree. Even so, the next theorem 
provides something similar but weaker. A system of sets is called \textbf{laminar} if any two members of it are either disjoint or 
$ \subseteq $-comparable.

\begin{thm}\label{GH lamináris van}
If $ b $ satisfies conditions 0,1,2,3',4, then there is a laminar system $ \mathcal{L}^{*} $ of optimal cuts such that any pair from $ V 
$ is separated optimally by some element of $ \mathcal{L}^{*} $.
\end{thm}
   
\begin{sbiz}
\begin{claim}
For any $ u\neq v\in V $ there exists an $ u-v $ cut $ X^{*} $ with $b(X^{*})= \lambda_b(u,v)$.
\end{claim}
 
\begin{ssbiz}
Let $ u,v $ be fix. The \textbf{error} of the sequence $ (X_n) $  of $ u-v $ cuts is  

\[\sum_{n=0}^{\infty} (b(X_n)-\lambda(u,v)).  \]

\noindent It is enough to prove the existence of a 
nested sequence $ (Y_n) $ of $ u-v $ cuts with finite error. Indeed, from the finiteness of the error it follows that $ 
\lim b(Y_n)=\lambda(u,v) $, hence by property 3'  \[\lambda_b(u,v)\leq b\left( \bigcap_{n=0}^{\infty}Y_n  \right)\leq \liminf 
b(Y_n) =\lim 
b(Y_n)=\lambda_b(u,v). \]

\begin{prop}\label{GH elsőfelfúj}
 For any sequence $ (X_n) $ with finite error there is another sequence $ (Z_n) $ with 
less or equal error such that $ Z_0 \supseteq \bigcup_{n=1}^{\infty}Z_n $.
\end{prop}  
\begin{sbiz}
Replace in the sequence $ (X_n) $ the 
member $ X_0 $ by $ X_0\cup X_1 $ and the member $ X_1 $ by $ X_1\cap X_0 $. By submodularity the 
error of the new sequence $ (X^{1}_n) $ is less or equal. Then replace  $ X^{1}_0=X_0\cup X_1 $  by $X^{2}_0:= X^{1}_0\cup X^{1}_2 = 
X_0\cup X_1\cup 
X_2 $ and replace $ X^{1}_2 $ by $X^{2}_2:=X^{2}_2\cap X^{2}_0= X_2\cap (X_0\cup X_1) $. In general let
 \[X^{m+1}_n= \begin{cases} X_{0}^{m}\cup X_{m+1}^{m} &\mbox{if } n= 0 \\
X_{0}^{m}\cap X_{m+1}^{m} & \mbox{if } n=m+1\\
X_n^{m} & \mbox{otherwise. } 
\end{cases}  \]

 Finally we claim that the following ``limit'' of these sequences is appropriate.  \begin{align*}
& Z_0:=\bigcup_{n=0}^{\infty}X_n\\
& Z_{n+1}:=X_{n+1}\cap \left( \bigcup_{i=0}^{n}X_n  \right).
\end{align*}
For 
\[ S_m:=\sum_{n=0}^{\infty} (b(X_n^{m})-\lambda(u,v)) \]
$ (S_m) $ is a non-negative decreasing  sequence thus it has a limit $ S $ i.e.
\[ S=\lim_{m\rightarrow \infty}\sum_{n=0}^{\infty} (b(X_n^{m})-\lambda(u,v)). \]
Consider the counting measure on $ \mathbb{N} $ and apply Fatou lemma:

 \begin{align*}
 S=&\liminf_{m}\sum_{n=0}^{\infty} (b(X_n^{m})-\lambda(u,v))\\
 \geq&\sum_{n=0}^{\infty} (\liminf_{m } b(X_n^{m})-\lambda(u,v))\\
 =&\liminf_m b\left( \bigcup_{i=0}^{m}{X_i}  \right)-\lambda(u,v)+ \sum_{n=1}^{\infty} (b(Z_n)-\lambda(u,v))\\
 \geq& \sum_{n=0}^{\infty} (b(Z_n)-\lambda(u,v)).
 \end{align*}
 
 At the last step we applied property 3' ($ S_0+\lambda(u,v) $ is an obvious upper bound for $ \{ b(X_0^{m}): m\in\mathbb{N} \} $). 
 Hence the 
 error of $ (Z_n) 
 $ is 
 smaller or equal 
 than the error of the earlier sequences.  
\end{sbiz}\\

Let $ X_n $ be a $ u-v $ cut with $ b(X_n)-\lambda(u,v)\leq 1/2^{n+1} $. Then the error of $ (X_n) $ is at most $ 1 $. Apply  
Proposition \ref{GH elsőfelfúj} with $ (X_n) $ to obtain $ (X_n^{1}) $ and let $ Y_0=X_0^{1} $.  Use Proposition \ref{GH elsőfelfúj} to 
the 
terminal segment of $ 
(X_n^{1}) $ consists all but the $ 0 $-th element (this sequence has error at most $ 1-(b(Y_0)-\lambda_b(u,v)) $)  to obtain $ (X_n^{2}) 
$ and let $ Y_1=X_0^{2} $. By continuing the process 
recursively  we build up a desired nested $ (Y_n) $ with error at most $ 1 $.
\end{ssbiz}

\begin{rem}
One can observe in the proof above that if $ (X_n) $ is a sequence of $ u-v $ cuts with finite error, then simply $ 
\bigcup_{n=0}^{\infty}X_n $ 
and $ 
\bigcap_{n=0}^{\infty}X_n  $ are optimal $ u-v $ cuts.
\end{rem}

 \begin{prop}\label{GH optvágásom metszete uniója opt vágás}
 The intersection and the union of (even infinitely many)  optimal $ u-v $ cuts is an optimal $ u-v $ cut.
 \end{prop}
 
 \begin{sbiz}
Let $ X $ and $ Y $ be optimal $ u-v $ cuts.  On the one hand,  $ b(X)\leq b(X\cup Y) $ and $ b(Y)\leq 
b(X\cap Y) $ hold since $ 
 X\cup Y $ and 
 $ X\cap Y $ are $ u-v $ cuts. Thus  \[ b(X)+b(Y)\leq b(X\cup Y)+b(X\cap Y). \]  On the other hand,  \[ b(X)+b(Y)\geq b(X\cup 
 Y)+b(X\cap Y) \]  by submodularity. Hence  equality holds and therefore  $ b(X)= 
 b(X\cup Y) $ and $ b(Y)= b(X\cap Y) $ since $ b(X), b(Y)<\infty $ because of property 4. By induction we know the statement for finitely 
 many 
 optimal $ u-v $ cuts. Consider an infinite 
 family $ \mathcal{X} $ of optimal $ 
 u-v $ cuts.  Let $ V=\{ v_n: n\in \mathbb{N} \} $ and let $ X_n' \in \mathcal{X} $ with $ v_n\notin X_n' $ if $ v_n\notin 
 \bigcap\mathcal{X}$  and an arbitrary element of $ \mathcal{X} $ otherwise.   
 Then $ X_n :=\bigcap_{m=0}^{n} X'_m $ is an optimal $ u-v $ cut again and  $  \bigcap_{n=0}^{\infty} 
 X_n =\bigcap \mathcal{X} $ as well by property 3'. 
 \end{sbiz}
 
 \begin{cor}\label{legkisebb optkut}
 There is a $ \subseteq $-smallest (largest) optimal $ u-v $ cut $ \boldsymbol{X_{u,v}}\ (Y_{u,v}) $ which is the intersection (union)  
 of 
 all optimal $ 
 u-v $ cuts.  
 \end{cor}

 \begin{claim}\label{esetekGH}
 Let  $ X $ be an optimal $ s-t $ cut and let $ Y $ be an optimal $ u-v $ cut. 
 \begin{enumerate}
 \item Assume  $ X $ is a $ u-v $ cut. Then $ Y\cup X $ is an optimal $ u-v $ cut if $ t\notin Y $ and $ Y\cap X $  is an optimal $ u-v $ 
 cut if $ t\in Y $.
 
 \item Assume $ X $ is a $ v-u $ cut. Then $  Y\cup(V\setminus X) $ is an optimal $ u-v $ cut if $ s\notin Y $ and $ Y \setminus X $ is 
 an 
 optimal $ u-v $ cut if $ s\in Y $. 
 \item Assume $ u,v\in X $. Then $ Y\cap X $ is an optimal $ u-v $ cut if $ t\notin Y $ and $ Y\cup(V\setminus X) $ is an optimal $ u-v $ 
 cut if $ t\in Y $.
 \item Assume $ u,v\notin X $. Then $ Y \setminus X  $ is an optimal $ u-v $ cut if $ s\notin Y $ and $  Y\cup X $ is an optimal $ u-v $ 
 cut if $ s\in Y $.
 \end{enumerate}
 \end{claim}

 \begin{sbiz}
 It is enough to prove 1 and 3 since  by replacing $ X $ with the optimal $ t-s $ cut $ V\setminus X $ in them we obtain 2 and 4 
 respectively. 
 To prove 1  assume first that 
 $ t\notin Y $. Since $ 
 X\cup Y $ is a $ s-t $ cut and $ X\cap Y $ is a $ u-v $ cut we have $ b(X\cup Y) \geq b(X)  $ and $ b(X\cap Y)\geq b(Y) $. 
 Combining this with  submodularity we get
 
 \[ b(X)+b(Y) \geq b(X\cup Y) + b(X\cap Y) \geq b(X)+ b(Y), \]
 
 \noindent thus  equality must hold in both  inequalities. If $ t\in Y $ and $ 
 s\in Y $, then $ X\cup Y $ is a $ u-v $ cut 
 and 
 $ X\cap Y $ is a $ s-t $ cut; therefore by arguing similarly as above we obtain that $ X\cup Y $ must be an optimal $ u-v $ cut. 
 Finally  
 if $ t\in Y $ and $ s\notin Y $, then on the one hand $ Y $ separates $ t $ and $ s $ and $ X $ does this optimally therefore $ b(X)\leq 
 b(Y) $, on the other hand $ Y $ is an optimal $ u-v $  cut and $ X $ is an $ u-v $ cut hence $ b(Y)\leq b(X) $. Thus $ 
 b(X)=b(Y) 
 $ therefore $ X $ and $ Y $ 
  both are optimal $ u-v $ cuts hence by Proposition \ref{GH optvágásom metszete uniója opt vágás}  $ X\cup Y $ and $ X\cap Y $ as well. 
 The proof of 3 is similar.
 \end{sbiz}

 \begin{cor}\label{belülszétválaszt}
 If $ X $ is an optimal $ s-t $ cut and $ u\neq v\in X $, then either $ X_{u,v} \subseteq X $ or $ X_{v,u}\subseteq X $ (where $ X_{x,y} 
 $ 
  stands for the $ \subseteq $-smallest optimal $ x-y $ cut). 
 \end{cor}
 
 \begin{sbiz}
  If $ X_{u,v} \subseteq X $, then  we are done. Assume $ X_{u,v} \not\subseteq X 
  $. By the minimality of $ X_{u,v} $ the $ u-v $ cut $ X_{u,v}\cap X $ cannot be optimal therefore by Claim 
   \ref{esetekGH}/3   $ X_{u,v}\cup(V\setminus X) $ is an optimal $ u-v $ cut. But then $ V \setminus [X_{u,v}\cup(V\setminus 
   X)]=X\setminus X_{u,v}  $ is an optimal $ v-u $ cut therefore we obtain $ X_{v,u} \subseteq (X \setminus X_{u,v})\subseteq X $.
 \end{sbiz}
 
 Theorem \ref{GH lamináris van} follows immediately from the next lemma (actually we need the lemma just with finite $ \mathcal{L} $).
 \begin{lem}\label{laminárisépítő}
 If $ \mathcal{L} $ is a  laminar system of optimal cuts and $ u\neq v\in V $, then there is a cut $ X^{*} $ for which $ 
 \mathcal{L}\cup \{ X^{*} \} $ is laminar and $ X^{*} $ separates optimally $ u $ and $ v $.  
 \end{lem} 
  \begin{sbiz}
 Let us partition $ \mathcal{L} $ into four parts $ \boldsymbol{\mathcal{L}_{u, \overline{v}}}:=\{ X\in \mathcal{L}:\ u\in X \wedge 
 v\notin X \} $, we 
 define $ \boldsymbol{\mathcal{L}_{\overline{u}, v}}, 
 \boldsymbol{\mathcal{L}_{u, v}} $ and $ \boldsymbol{\mathcal{L}_{\overline{u}, \overline{v}}} $ similarly. 
 If $ X_{u,v} \subseteq \widehat{X} $  for some $ \widehat{X}\in 
 \mathcal{L}_{u,\overline{v}} $, 
 then  $ \{X_{u,v} \}\cup \mathcal{L}_{u,v}\cup \mathcal{L}_{\overline{u},v} $ is laminar. Suppose that we have no such an $ \widehat{X} 
 $ not  even if we 
 interchange $ u $ and $ v $. By Corollary \ref{belülszétválaszt} we know that for all $ W\in 
 \mathcal{L}_{u,v} $ 
 either $ X_{u,v} \subseteq W 
 $ or $ X_{v,u} \subseteq W $. Hence by symmetry we may assume that $   \bigcap \mathcal{L}_{u,v}\supseteq X_{u,v} $. We will show 
 that 
  $ \{ X_{u,v} \}\cup \mathcal{L}_{u,v}\cup \mathcal{L}_{\overline{u},v} $ is laminar in this case as well. Let $ X\in 
  \mathcal{L}_{\overline{u}, v} $ 
 be arbitrary. Then $ X_{v,u} \not\subseteq X $ otherwise $ \widehat{X}:=X $ would be a bound. But then $ X_{v,u}\cap X $ cannot be an 
 optimal $  v-u  $ cut by the minimality of $ X_{v,u} $. Therefore by Claim \ref{esetekGH}/1 we know that $ X_{v,u}\cup X $ is an optimal 
 $ v-u $ cut 
 and hence $ 
 V\setminus (X_{v,u}\cup X) $ is an optimal $ u-v $ cut. Thus $ V\setminus (X_{v,u}\cup X) \supseteq  X_{u,v} $ from which $ X\cap 
 X_{u,v}=\varnothing $ follows.
 
  Thus we may suppose that $ \{ X_{u,v} \}\cup \mathcal{L}_{u,v}\cup \mathcal{L}_{\overline{u},v} $ is laminar. If for some 
  $ Y\in \mathcal{L}_{u,\overline{v}} $ the set $ \{ X_{u,v}, Y \} $ is not laminar, then the cut
  $ X_{u,v}\cap Y $ may not be an optimal $ u-v $ cut because of the definition of $ X_{u,v} $. But then $ X_{u,v}\cup Y $ is an optimal 
  $ u-v $ cut by Claim 
  \ref{esetekGH}/1. Let  \[\mathcal{Y}:= \{ Y: Y\in \mathcal{L}_{u,\overline{v}}\wedge \{ X_{u,v}, Y \}\text{ is not 
  laminar}\}.  \] 
  The set $ \{ X_{u,v}\cup Y: Y\in \mathcal{Y} \} $ consists of optimal $ u-v $ cuts and  totally ordered by $ \subseteq $. By taking a 
  cofinal sequence of type at most 
  $ \omega $ and applying 3'  we obtain that $X_0:= \bigcup \mathcal{Y} $ is an optimal $ u-v $ cut. Note that 
  $ \{ X_0 \}\cup( \mathcal{L}\setminus \mathcal{L}_{\overline{u},\overline{v}}) $ is laminar. For each 
  $ Z\in \mathcal{L}_{\overline{u}, \overline{v}} $ fix some $ s_Z,t_Z $ such that $ Z $ is an optimal $ s_Z-t_Z $ cut. We claim that if 
  for such a 
  $ Z $ the pair $ \{ X_0,Z \} $ is not laminar, then $ \{ X_{u,v},Z \} $ is not laminar as well. Indeed, 
  consider just the construction of $ X_0 $ and the fact that if for a cut $ Y\in \mathcal{L}_{u,\overline{v}} $ we have $ Y\cap Z\neq 
  \varnothing $, then $ Z\subseteq Y $ by the laminarity. Let
  \[ \mathcal{Z}:=\{ Z\in \mathcal{L}_{\overline{u}, \overline{v}}: \{X_0, Z \}\text{ is not laminar} \}. \]
  For $ Z\in \mathcal{Z} $ we know that $ \{ X_{u,v},Z \} $ is not laminar. By the 
  definition of $ X_{u,v} $ the cut $ 
  X_{u,v}\setminus Z $ may not be an optimal $ u-v $ cut hence by Claim \ref{esetekGH}/4 it follows, that $ s_Z\in X_{u,v}(\subseteq X_0) 
  $. 
  Finally by 
  Claim \ref{esetekGH}/4 we may take $X^{*}= X\cup \bigcup \mathcal{Z} $ by adding countably many elements of $ \mathcal{Z} $ with union 
  $ \bigcup \mathcal{Z} $  one by one to $ X $ and taking limit. 
  \end{sbiz}
\end{sbiz}

   \section{A counterexample}
   
   In the previous section we obtained (as a special case of Theorem \ref{GH lamináris van})  the existence of a laminar 
   system of optimal cuts for 
   countably infinite finitely separable 
   weighted graphs which elements separate any vertex pair optimally.  In this section we show by an example that one cannot 
   guarantee the existence of a Gomory-Hu tree as well without further assumptions.
   Let $ V= \{ v_n: n\in \mathbb{N} \} $ and let $ E=\{ v_{\infty}v_n:n\in \mathbb{N}  \}\cup \{ v_nv_{n+1}: n\in \mathbb{N} \} $. 
   Finally $ 
   c(v_{\infty}v_n):=1 $ 
   for all 
   $ n\in \mathbb{N} $ and with the notation $ e_n:=v_nv_{n+1} $

   \[ c(e_n):= 
   \begin{cases} 2 &\mbox{if } n=0 \\
   c(e_{n-1})+n+1 & \mbox{if } n>0.  
   \end{cases} \]

   \begin{figure}[H]
   \centering
   
   \begin{tikzpicture}

   \node[circle,inner sep=0pt,draw,minimum size=5] (v1) at (0,0) {$v_0$};
    \node[circle,inner sep=0pt,draw,minimum size=5] (v2) at (2.5,0) {$v_1$};
    \node[circle,inner sep=0pt,draw,minimum size=5] (v3) at (5,0) {$v_2$};
    \node[circle,inner sep=0pt,draw,minimum size=5] (v4) at (7.5,0) {$v_3$};
    \node (v6) at (10,0) {$\dots$};
    \node[circle,inner sep=0pt,draw,minimum size=5] (v5) at (5,2.5) {$v_\infty$};
    \draw  (v1) edge (v2);
    \draw  (v2) edge (v3);
    \draw  (v3) edge (v4);

    \draw  (v1) edge (v5);
    \draw  (v2) edge (v5);
    \draw  (v3) edge (v5);
    \draw  (v4) edge (v5);
   \draw  (v4) edge (v6);
   
   \draw  (v5) edge (v6);

   \node at (2,1.2) {$1$};
   \node at (3.4,1.2) {$1$};
   \node at (4.8,1.2) {$1$};
   \node at (6.6,1.2) {$1$};
   \node at (8,1.2) {$1$};

   \node at (1.2,-0.2) {$2$};
   \node at (3.8,-0.2) {$4$};
   \node at (6.2,-0.2) {$7$};
   \node at (8.8,-0.2) {$11$};

   \end{tikzpicture}
   \caption{ A finitely separable weighted graph without a Gomory-Hu tree} 
   
   \end{figure}
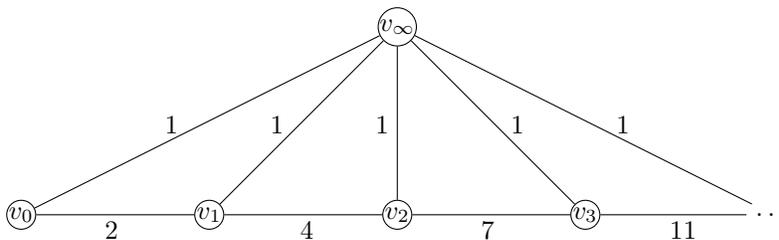

  \begin{claim}
  If $ n<m $, then $ \{ v_0, v_1,\dots,v_n \}=:V_n $ is the only optimal $ v_n-v_m $ cut.
  \end{claim}
  
  \begin{sbiz}
  Pick an optimal $ v_n-v_{m} $ cut $ X $. Since $ d_c(V_n)<c(e_k) $ whenever $ k>n $, a cut $ X $ may not 
  separate the end vertices of such an $ e_k $. Then $ v_\infty\notin X $ otherwise $ d_c(X)=\infty $. Thus we have $ X \subseteq V_n $. 
  Suppose, for seeking a contradiction, that $ v_l\notin X $ for some $ l<n $ and $ l $ is the largest such an  index. Then
  
  \[ d_c(X) -d_c(V_n) \geq c(e_l)-l-1>0, \]
  which contradicts  to the optimality of $ X $.
  \end{sbiz}
  
  \begin{claim}
  $ (G,c) $ has no Gomory-Hu tree.
  \end{claim}
   \begin{sbiz}
   Assume, to the contrary, that $ T $ is a Gomory-Hu tree of $ (G,c) $. For all $ e\in E(T) $ pick the  fundamental cut $ X_e $ that 
   corresponds 
   to 
   $ e $ and does not contain $ v_\infty $. On the one hand, $ \mathcal{L}:=\{ X_e \}_{e\in E(T)} $ is a laminar system of optimal cuts 
   that 
   contains at least one 
   $ \subseteq $-maximal element (if $ e $ incident with $ v_\infty $ in $ T $, then $ X_e $ is a maximal element). On the other hand, $ 
   \mathcal{L}=\{ V_n: n\in \mathbb{N} \} $ since the optimal cuts are unique up to complementation and the additional condition ``does 
   not contain $ v_{\infty} $'' 
   makes 
   them 
   unique. This is a contradiction since $ (V_n) $ is a strictly increasing sequence.
   \end{sbiz}
   
  \begin{rem}
   One can obtain also a locally finite counterexample by some easy modification of our counterexample above. 
  \end{rem}
 
 \section{Existence of an abstract Gomory-Hu tree}

In this section we prove or main result (which is Theorem \ref{GH van absztract GHtree}).
It will be convenient to use the following  equivalent but formally weaker definition of Gomory-Hu trees.
    \begin{claim}\label{GH fa def gyengít}
    $ T=(V,F) $ is a Gomory-Hu tree with respect to $ b $ if for all $ uv\in F $ the fundamental cuts corresponding to $ uv $ in $ T $ 
    separate 
    optimally $ u 
    $ and $ v $.
    \end{claim}
   \begin{sbiz}
    Let $ u\neq v\in V $ be arbitrary and let $ v_1,v_2,\dots, v_m $ be the vertices of the unique 
    $ u-v $ path in $ T $ numbered  in the path order.
    
    \begin{prop}
     For all pairwise distinct $  u,v,w\in V$ we have: \[ \lambda_b(u,w) \geq \min\{ \lambda_b(u,v), \lambda_b(v,w) \}. \]
    \end{prop}
    
    \begin{nbiz}
    It follows from the fact that if a cut separates $ u $ and $ w $, 
    then it separates either $ u $ and $ v $ or $ v $ and $ w $ as well. \rule{1ex}{1ex}
    \end{nbiz}\\
    
     On one hand by applying the Proposition above repeatedly we obtain
    
      \[ \lambda_b(u,v)\geq 
    \min\{ \lambda_b(v_i,v_{i+1}):\ 1\leq i<m \}=:\lambda_b(v_{i_0},v_{i_0+1})\text{ for some } 1\leq i_0 <m. \]   
    
    On the other hand, the fundamental cuts corresponding to the edge $ v_{i_0}v_{i_0+1} $  
    separates $ u $ and $ v $ and have value $ \lambda_b(v_{i_0},v_{i_0+1}) $ by assumption. Thus
      
      \[ \lambda_b(u,v)\leq\lambda_b(v_{i_0},v_{i_0+1}). \]

    Hence equality holds and
    the fundamental cuts corresponding to  $ v_{i_0}v_{i_0+1}\in F $ are optimal cuts between $ u $ and $ v $.
   \end{sbiz}\\

 A sequence $ (X_n) $ of optimal cuts is \textbf{essential} if all of its members separate optimally a vertex pair that the 
 earlier members do not. 
 \begin{lem}\label{monoton vágás sorozatok}
 If $ (X_n) $ is a $ \subseteq $-monotone sequence of optimal cuts and $\lim X_n =: X \notin \{ \varnothing, V \} 
 $, then $(X_n) $ has no essential subsequence.    
 \end{lem}
 
 \begin{ssbiz}
 Assume, to the contrary, that  $ (X_n) $ is a counterexample. By symmetry we may suppose that $ (X_n) $ is increasing. By trimming $ 
 (X_n) $ we may assume that it is essential witnessed by $ s_n,t_n $ i.e. $ X_n $ is an optimal $ s_n-t_n $ cut but $ X_m $ is 
 not 
 whenever $ m<n $.

 \begin{claim} \label{t_n künn van}
  $ t_n\notin X $ holds for all large enough $  n $.
 \end{claim}
 
 \begin{sbiz}
Suppose that $ t_n\in X $ for infinitely many $ n $. By the monotone-continuity the numerical sequence $ 
(b(X_n)) $ converges to $ 
b(X)>0 $, thus $ b(X_n)\geq b(X)/2>0 $ for large 
enough $ n $. On the one hand,
  $ b(X \setminus X_n)\rightarrow 0 $ since $ (X\setminus X_n) \rightarrow \varnothing $ monotonously. On the other hand, $ X 
  \setminus 
  X_n $ is a $ t_n-s_n $ cut for infinitely many $ n $ because of the indirect assumption and for such an $ n $
   \[ b(X \setminus X_n)\geq b(X_n)\geq b(X)/2>0 \]
   which is a contradiction.   
 \end{sbiz}\\
 
  By trimming $ (X_n) $, we may assume that $ t_n\notin X $ for all $ n $. It implies that $ b(X_n)\leq b(X_{n+1}) $ for 
  each $ n $ because $ X_{n+1} $ is an $ s_n-t_n $ cut and $ X_n $ is an optimal $ s_n-t_n $ cut. But then $ s_{n+1}\notin 
  X_n 
  $  for all $ n $, otherwise $ X_n $ would be at least as good $  s_{n+1}-t_{t_{n+1}} $ cut as the optimal one but $  X_{n+1} $ 
  is the first optimal $ s_{n+1}-t_{n+1} $ cut of the sequence by the choice of $ (X_n) $.

  \begin{claim}
  $ X_n $ is an optimal $ s_n-s_{n+1} $ cut for all $ n $.
  \end{claim}
  
  \begin{sbiz}
  By Corollary \ref{belülszétválaszt}, there is an $ Y\in \{ X_{s_n,s_{n+1}}, X_{s_{n+1},s_n} \} $ for which $ Y\subseteq X_{n+1} $. If $ 
  b(Y)<b(X_n) 
  $ would hold, 
  then $ Y $ would be either a better $ s_n-t_n $ cut than $ X_n $ or a better $ s_{n+1}-t_{n+1} $ cut than $ X_{n+1} $; which both are 
  impossible.
  \end{sbiz}\\
  
  \noindent Since $ s_n\in X $ for al $ n $ the Claim above contradicts  to Claim \ref{t_n künn van} with the choices $ 
  s_n:=s_n $ and $ 
  t_n:=  s_{n+1} $.
 \end{ssbiz}\\

  Take an optimal cut $ X $. For $ u\neq v\in X $ let $ \boldsymbol{u\prec_X v} $ if $ 
  X_{u,v}\not  \subseteq X $.

  \begin{claim}
  The relation $ \prec_X $ is a strict partial ordering on $ X $.
  \end{claim}
  
  \begin{sbiz}
  It is irreflexive  by definition. For transitivity assume $ u\prec_X v\prec_X w $. If $ u=w $, then we have  $ u\prec_X v $ and $ 
  v\prec_X u $ 
  which 
  contradicts to 
   Corollary \ref{belülszétválaszt}. Thus we may assume that $ u,v,w $ are pairwise distinct. Suppose, to contrary, that $ u\prec_X w 
  $ 
  does not hold i.e.  $ X_{u,w}\subseteq X $. Assume first that 
   $v\in X_{u,w} $. By Corollary \ref{belülszétválaszt}, either $ X_{u,v} \subseteq X_{u,w} $ or $ X_{v,u}\subseteq X_{u,w} $.  Since $ 
   u\prec_X v $, necessarily $ X_{v,u}\subseteq X_{u,w} $. But then $ X_{u,w} $ and $ X_{v,u} $ are both $ v-w $ cuts and
  
  \[ \lambda_b(v,w) \geq \min \{ \lambda_b(v,u),\ \lambda_b(u,w) \}=\min \{ b(X_{v,u}),\ b(X_{u,w}) \}  \]
  
  \noindent shows that one of them is optimal  which contradicts to  $ v\prec_X w $.
  
   Hence $ v\notin X_{u,w} $ holds. $ X_{u,w} $ is not an optimal $ u-v $ cut since $ u\prec_X v $. Therefore  $ 
   b(X_{u,w})>b(X_{v,u}) 
   $. 
 (Note that $ X_{v,u}\subseteq X $ by $ u\prec_X v $  and  by Corollary \ref{belülszétválaszt}). Hence $ w\notin X_{v,u} $ otherwise $ 
 X_{v,u} $ 
  would be  a 
  better cut between 
  $ w $ and $ u $ than the optimal. On the other hand, $ X_{v,u} $ is not an optimal $ v-w $ cut since $ v\prec_X w $ hence  $ X_{w,v} 
  \subseteq 
  X $ and $ b(X_{w,v})<b(X_{v,u}) $ hold. Necessarily $ u\in X_{w,v} $, otherwise $ X_{w,v} $ separates better $ w $ and $ u $ than $ 
  X_{u,w} $, but 
  then $ X_{w,v} $ separates better $ u $ and $ v $ than $ X_{v,u} $, which is a contradiction.
  \end{sbiz}

  \begin{lem}\label{minimális pont ralizál}
  If $ X $ is an optimal $ s-t $ cut, then $ X $ has a $ \prec_X  $-minimal element $ s' $. For all such an $ s' $,  cut $ X $ it is 
  an 
  optimal 
  $ s'-t $ cut.
  \end{lem}

  \begin{ssbiz}
  Let $ A=\{ x\in X: \lambda_b(x,t)=  \lambda_b(s,t) \} $ and $ B:= \{ y\in X: \lambda_b(y,t)<  \lambda_b(s,t) \}$. Then $ 
  A\cup B  $ is a partition of $ X $. Note that $ A\neq \varnothing $ since $ s\in A $.

   \begin{prop}\label{a<b}
   For all  $ x\in A $ and $ y\in B:\  x\prec_X y $ holds.
   \end{prop}
   
   \begin{sbiz}
   If $ x\in A $ and $ y\in B $, then $ \lambda_b(x,y) <\lambda_b(s,t)(=\lambda_b(x,t))$, otherwise  
   
   \[ \lambda_b(y,t)\geq \min \{ \lambda_b(x,y), \lambda_b(x,t) \}=\lambda_b(x,t)=\lambda_b(s,t)  \]

 \noindent contradicts to $ y\in B $. Therefore if $ X_{x,y} \subseteq X $ would hold, then (since  $ X_{x,y} $ is a $ x-t $ cut)
 \[ \lambda_b(x,t)  \leq \lambda_b(x,y) <\lambda_b(s,t)  \]
 
 \noindent which is impossible since $ x\in A $. 
 \end{sbiz}\\

   By  Proposition \ref{a<b}, it is enough to find a minimal element for the poset $ (A,\prec_X ) $. The existence of such an element 
   follows 
   immediately from the following Proposition.

   \begin{prop}\label{B véges}
   Set $ A $ is finite. 
   \end{prop}
   
   \begin{sbiz}
    Assume, to seeking for contradiction, that $ A $ is infinite. Pick a nested sequence $ (A_n) $ of nonempty subsets of $ A $ with $ 
    \bigcap_{n=0}^{\infty}A_n=\varnothing $. On the one hand, $ b(A_n) \rightarrow 0$ by property 3. On the other hand, 
    every $ A_n $ separates an $x\in A $ from $ t $ and hence \[ b(A_n) \geq\lambda_b(x,t)=\lambda_b(s,t)>0, \]
    which is a contradiction.
   \end{sbiz}\\
   
   For the second part of  Lemma \ref{minimális pont ralizál} let $ s' $ be a $ \prec_X  $-minimal element of $ X $.  Then by Proposition 
   \ref{a<b} $ s'\in A $ thus $ \lambda_b(s',t)=\lambda_b(s,t) $ by the definition of $ A $.
  \end{ssbiz}

 \begin{claim}
 For any $ s\in V $ the family  $\boldsymbol{\mathcal{C}_s}:= \{ X_{u,s}: u\in V\setminus \{ s  \} \} $ of optimal cuts is laminar.
 \end{claim} 
 \begin{sbiz}
 Suppose, to the contrary, that $ \{ X_{u,s},\ X_{v,s} \}\subseteq \mathcal{C}_s $ is not laminar.  If $ u\in X_{v,s} $, then $ 
 X_{u,s}\cap X_{v,s} $ is an $ u-s $ cut and $ X_{u,s}\cup X_{v,s} $ is a $ v-s $ cut. By 
 submodularity $ X_{u,s}\cap X_{v,s} $ is an optimal $ u-s $ cut (and $ X_{u,s}\cup X_{v,s} $ is an optimal $ v-s $ cut) which 
 contradicts to the definition of $ X_{u,s} $. For  $ v\in X_{u,s} $ the argument is the same. Finally if $ u\in 
 X_{u,s}\setminus X_{v,s} $ and $ v\in X_{v,s}\setminus X_{u,s} $, then $ X_{u,s}\setminus X_{v,s} $ is an $ u-s $ cut and $ 
 X_{v,s}\setminus X_{u,s} $ is an $ v-s $ cut thus by applying Remark \ref{GH szubmodegye 2} follows that they are also optimal, 
 contradicting 
 to the definitions of $ X_{v,s}$ and $ X_{u,s} $.
 \end{sbiz}\\
 
 Let $ \prec_V $ be the trivial partial ordering on  $ V $ (i.e. under which there are no comparable elements).
 \begin{lem}\label{particionál levél}
 Let $ X $ be either an optimal  cut or $ V $. Pick an $ \prec_X  $-minimal  element $ s $ of $ X $ (see Lemma \ref{minimális pont 
 ralizál}). Then the $ 
 \subseteq $-maximal elements of 
 the laminar system  
  $\boldsymbol{\mathcal{C}_{s,X}}:= \{ X_{u,s}: u\in X\setminus \{ s  \} \} $  forms a partition of $ X\setminus \{ s \} $. 
 \end{lem}
 \begin{sbiz}
  By the choice of $ s $ we know that $ \bigcup \mathcal{C}_{s,X}\subseteq X\setminus \{ s \} $. It is enough to show that  for any $ 
  u^{*}\in X\setminus \{ s \} $ the laminar system  
 $\mathcal{C}_{s,X} $  has a maximal element that contains $ u^{*} $. Assume, seeking for 
 contradiction, that it is false and $  (X_{u_n,s}) $ is a strictly increasing sequence that shows this. On the one hand, this sequence 
 is essential 
 because $ X_{u_m,s} $ may not be an optimal $ u_n-s $ cut for $ m<n $ since $ X_{u_n,s} $  is the $ \subseteq $-smallest such a cut. On 
 the other 
 hand, $ \lim X_{u_n,s} \subseteq V\setminus \{ s \}$ which contradicts to Lemma \ref{monoton vágás sorozatok}.
 \end{sbiz}\\

 We build  the desired abstract Gomory-Hu tree for $ b $ by using  Lemma \ref{particionál levél}  repeatedly. Pick an arbitrary $ r\in V 
 $ for  root. It 
 makes possible to define a unique fundamental cut for each edge $ e $ of the tree, namely the vertex set of the component after deletion 
 of 
 $ 
 e $ that does not contain $ r $. Let $ \{ X_i \}_{i\in I_0} $ consists of the maximal elements of the laminar system 
 $ \mathcal{C}_{r} $. Let $ x_i $ be a $ \prec_{X_i} $-minimal element of $ X_i $ and draw the tree-edges $ rx_i $ for $ i\in I_0 $. 
 Note that Lemma \ref{minimális pont ralizál} ensures that the fundamental cut corresponds to $ rx_i $ will separate optimally $ r $ and 
 $ x_i 
 $ assuming that $ X_i $ will be the vertex set of the subtree rooted at $ x_i $. 
 Take now for each $ i\in I_{0} $ the  $ \subseteq $-maximal elements $ \{ X_{i,j} \}_{j\in I_1} $ of  
 $ \mathcal{C}_{x_i,X_i} $ and choose a $ \prec_{X_{i,j}} $-minimal element $ x_{i,j} $ of $ X_{i,j} $. Draw the 
 tree-edges $ x_ix_{i,j} $ for all $ i\in I_0 $ and $ j\in I_1 $. By continuing recursively the process we claim that every $ v\in V $ 
 has to appear in 
 the 
 tree. Indeed, if some $ v $ does not, then we would obtain a nested essential sequence of optimal cuts such that its limit 
 contains  $ v $ which contradicts to Lemma \ref{monoton vágás sorozatok}.

\end{document}